\documentclass[a4paper, svgnames, 11pt]{amsart}
\usepackage[english]{babel}
\usepackage{amsmath, amsthm, amssymb, hyperref, color}
\usepackage[capitalize,noabbrev]{cleveref}
\usepackage{graphicx}
\usepackage[dvipsnames]{xcolor}
\usepackage{caption}
\usepackage{subcaption}
\usepackage{mathtools}
\usepackage{enumerate}
\usepackage{verbatim}
\usepackage{tikz}
\usepackage{tikz-cd}
\usepackage[inline,shortlabels]{enumitem}
\usepackage{hyperref}
\hypersetup{
 colorlinks=true,
 linkcolor=blue,
 filecolor=blue, 
 urlcolor=blue,
 citecolor=blue,
}
\usepackage{url}
\usepackage{longtable}
\usepackage{bm}
\usepackage{xypic}
\usepackage{soul}
\usepackage{stmaryrd}
\usepackage{subcaption}
\usepackage{hhline}
\usepackage[a4paper]{geometry}
\numberwithin{equation}{section}

\newtheorem{definition}{Definition}[section]
\newtheorem{lemma}[definition]{Lemma}
\newtheorem{conjecture}[definition]{Conjecture}
\newtheorem{theorem}[definition]{Theorem}
\newtheorem{corollary}[definition]{Corollary}

\newcommand{\cd}{\textnormal{cd}}

\title[Exceptional Groups]
{Exceptional groups of order $p^6$ for primes $p\geq 5$} 

\author{E.A.\ O'Brien}
\address{Department of Mathematics, University of Auckland, Auckland, New Zealand}
\email{e.obrien@auckland.ac.nz}
\author{Sunil Kumar Prajapati}
\address{Indian Institute of Technology Bhubaneswar, Arugul Campus, Jatni, 
 Khurda-752050, India}
\email{skprajapati@iitbbs.ac.in}
\author{Ayush Udeep$^*$}
\address{Indian Institute of Science Education and Research Mohali, Sector 81, SAS Nagar, Punjab 140306, India}
\email{udeepayush@gmail.com}
\thanks{$^{\textbf{*}}$Corresponding author.
}

\thanks{
O’Brien was supported by the Marsden Fund of New Zealand Grant 23-UOA-080
and by a Research Award of the Alexander von Humboldt Foundation. Prajapati acknowledges
the Science and Engineering Research Board, Government of India, for financial support through
grant MTR/2019/000118. Udeep thanks the Indian Institute of Science Education 
and Research Mohali for his Postdoctoral Fellowship.
}

\begin{document}
 
\maketitle 
 
\begin{abstract}
The minimal faithful permutation degree $\mu(G)$ of a finite group
$G$ is the least integer $n$ such that $G$ is isomorphic to
a subgroup of the symmetric group $S_n$.
If $G$ has a normal subgroup $N$ such that $\mu(G/N) > \mu(G)$, 
then $G$ is exceptional.
We prove that the proportion of exceptional groups of order $p^6$ for primes
$p \geq 5$ is asymptotically 0. We identify $(11p+107)/2$ such groups 
and conjecture that there are no others.
 \end{abstract}
 
 \section{Introduction}
The minimal faithful permutation degree $\mu(G)$ of a finite group
$G$ is the least integer $n$ such that $G$ is isomorphic to
a subgroup of the symmetric group $S_n$.
This invariant has attracted significant attention:
see, for example, \cite{EH, EP, EST, DLJ, karpilovsky1970, W}. 
Cannon, Holt, and Unger \cite{CHU19} exploit 
such representations to assist with structural
investigations of large degree matrix groups. 
 
Clearly, $\mu(H) \leq \mu(G)$ for every subgroup $H$ of $G$, 
but this need not hold for quotients of $G$. 
If $G$ has a normal subgroup $N$ such that $\mu(G/N) > \mu(G)$, then 
$N$ is a \emph{distinguished subgroup}, 
$G/N$ is a \emph{distinguished quotient}, and
$G$ is \emph{exceptional}. 
Easdown and Praeger \cite{EP} introduced this terminology, 
believing that such groups are rare, 
and showed that the smallest exceptional group has order 32. 
Kov\'{a}cs and Praeger \cite{KP00} 
proved  that $\mu(G/N) \leq \mu(G)$ if $G/N$
has no nontrivial abelian normal subgroup. Holt and Walton
\cite{HW02} proved that there exists a constant $c$ such that
$\mu(G/N) \leq c^{\mu(G) - 1}$ for every $N \lhd G$. 
Franchi \cite{F11} proved that if
$G$ is a finite non-abelian $p$-group with abelian maximal subgroup, then
$\mu(G/G') \leq \mu(G)$. 

Lemieux \cite{LemieuxExc} (and later
Chamberlain \cite{chamberlain2018minimal}) 
constructed an exceptional group of order $p^5$ for every odd prime $p$ 
and proved that there are no exceptional groups of order dividing $p^4$. 
Britnell, Saunders, and Skyner \cite{BSSexc} determined the $p + 6$ 
exceptional groups of order $p^5$ for $p \geq 5$; 
since the number of groups of order $p^5$ is at most $2p + 71$, 
the proportion of exceptional groups of order $p^5$ is asymptotically $1/2$. 
Hence they questioned whether ``exceptional" is an apposite
descriptor for such groups.  By contrast, we prove the following.
\begin{theorem}\label{thm:main}
The proportion of exceptional groups of order $p^6$ for primes
$p \geq 5$ is asymptotically $0$.  
\end{theorem}
We prove this theorem by obtaining in 
Theorem \ref{thm:upperboundexceptional} an 
upper bound to the number of such groups.
We identify $(11p+107)/2$ exceptional groups of order $p^6$ 
and conjecture that there are no others.

For completeness, we report that 32 of the 267 
groups of order $2^6$ are exceptional, as are 67 
of the 504 groups of order $3^6$.  This computation takes about 
10 minutes of CPU time using the {\sc SmallGroups} library 
\cite{SmallGroups} and {\sc Magma} 2.28-13 \cite{MAGMA} on a 2.6 GHz machine. 
The identifiers of these groups are recorded at \cite{github}.

\section{Background and notation} \label{sec:prelims}
Let $p\geq 5$ be prime. 
Let $\omega$ be the smallest positive integer which is 
a primitive root modulo $p$, and let $\nu$ be the smallest 
positive integer which is a quadratic non-residue modulo $p$.

Let $Z(G)$, $G'$ and $d(G)$ denote respectively 
the center, commutator subgroup and
 minimal number of generators of a $p$-group $G$. 
Let $\cd(G)$ and $\exp(G)$ denote respectively the set of 
irreducible character degrees and exponent of $G$. 

A {\it quasi-permutation matrix} is a square matrix over the complex field
with non-negative integral trace.
Wong \cite{W} defined a
{\it quasi-permutation representation} of $G$ as
one consisting of quasi-permutation matrices.
Let $c(G)$ denote the minimal degree of
a faithful quasi-permutation representation of $G$.
It is easy to deduce that $c(G)\leq \mu(G)$;
Behravesh and Ghaffarzadeh
\cite[Theorem 3.2]{BG} proved that
$c(G) = \mu(G)$ for a $p$-group $G$ of odd order.

Newman, O'Brien and Vaughan-Lee \cite{NO'bV2004} established that there are
\begin{equation}\label{Total}
3p^2 + 39 p + 344 + 24\gcd(p-1,3) + 11 \gcd(p-1,4) + 2 \gcd(p-1,5) 
\end{equation}
groups of order $p^{6}$ where $p \geq 5$.
The groups are classified into 43 isoclinism
families labelled $\Phi_{i}$ for $1 \leq i \leq 43$;
the first 10 also contain groups of order dividing $p^5$. 
Invariants for these families are tabulated in \cite{RJ};
these include the structures of the derived group
and the central quotient, and the character degrees.
Parameterized presentations are listed in \cite{Arxiv}
for the groups of order $p^6$ and in \cite{RJ} for orders dividing $p^5$.
(See \cite[Remark 3.3]{ESAp6} for those few presentation
of 5-groups not covered by \cite{Arxiv}.)
We use the identifiers of \cite{Arxiv} for groups 
of order $p^6$ and those of \cite{RJ} for quotients. 
The values of $c(G)$ and $\mu(G)$ for the groups $G$ of 
order $p^5$ and~$p^6$ are available  
in \cite{PUp5} and \cite{ESAp6} respectively. 

As is well known (see, for example, \cite[Theorem 2]{DLJ}),
if $G \cong C_{p^{n_1}} \times \ldots \times C_{p^{n_k}}$, then
$\mu(G) = \sum_{i=1}^k p^{n_i}$.  Clearly, no abelian group
is exceptional. 
All 11 abelian groups of order $p^6$ are in $\Phi_1$. 

\section{An upper bound to the number of exceptional groups 
of order $p^6$} \label{sec:upperbound}
We obtain an upper bound to the number of
exceptional groups of order $p^6$, and so prove \cref{thm:main}.

Recall from \cite{lewis2009character} that a group $G$ is {\it VZ} if $\chi(g)=0$ for all
$g\in G\setminus Z(G)$ and all non-linear characters $\chi$ of $G$.
\begin{lemma} \label{lemma:VZexc}
 Let $G$ be a VZ group of order $p^6$. If $d(Z(G)) = d(G')$ and 
$|Z(G)| = p^2$, then $G$ is not exceptional.
\end{lemma}
\begin{proof}
 From \cite[Corollary 4]{SA}, we deduce that 
 \begin{equation*}
 c(G) = 
 \begin{cases}
 2p^3 &\text{ if } Z(G) \cong C_p \times C_p\\
 p^4 &\text{ if } Z(G) \cong C_{p^2}.
 \end{cases}
 \end{equation*}
If $G$ has a quotient $Q$ of order $p^4$, then $Q$ is not distinguished. 
Let $Q = G/N$ where $|N| = p$.
Since $N \leq Z(G)$, 
 \[ \frac{G}{Z(G)} \cong \frac{G/N}{Z(G)/N}. \]
From \cite[Lemma 2.4]{lewis2009charactertable}, $G/Z(G)$ is elementary abelian. 
Since $Z(G)/N \cong C_{p}$, we conclude that 
$\exp(G/N) \leq p^2$. If $\cd(Q) = \{ 1, p \}$, then
$c(Q) \leq 2p^3$ (see \cite[Theorem 3]{SAcyclic}). 
If $\cd(Q) = \{ 1, p^2 \}$, then $Q\in \Phi_{5}$ (see
\cite[\S 4.1]{RJ}). Hence $Q$ is a VZ group
with center of order $p$, and $c(Q) = p^3$ (see \cite[Corollary 4]{SA}). 
If $\cd(Q) = \{ 1, p, p^2 \}$, then
 $Q\in \Phi_{i}$ where $i \in \{ 7, 8, 10 \}$ and $c(Q) = p^3$
(see \cite[\S 4.1]{RJ}).
\end{proof}

\begin{corollary} \label{cor:phi15}
No group in $\Phi_{15}$ is exceptional.
\end{corollary}
\begin{proof}
If $G \in \Phi_{15}$, then 
$d(Z(G)) = d(G')$ and $|Z(G)| = p^2$ (see \cite{Arxiv}). 
Lemma \ref{lemma:VZexc} shows that $G$ is not exceptional. 
\end{proof}

\begin{lemma}\label{lemma:phi21exc}
No group in $\Phi_{21}$ is exceptional.
\end{lemma}
\begin{proof}
If $G \in \Phi_{21}$, then 
$\exp(G) \leq p^2$ (see \cite[Table~1]{Arxiv}). 
By \cite[Theorem 3.2]{BG} and
\cite[Lemma 4.28]{ESAp6}, we deduce that $\mu(G) = 2p^3$. 
Let $1\neq N \lhd G$ and let $Q = G/N$.
If $Q$ is cyclic, then it is not distinguished 
(see \cite[Lemma 1.2]{BSSexc}). 
Suppose $Q$ is not cyclic. If $|Q| = p^4$, 
then $\mu(Q) \leq 2p^3$ (see \cite[Table 1]{LemieuxExc}. 
Suppose $|Q| = p^5$. If $Q$ is abelian,
then $\mu(Q) < 2p^3$ since $\exp(Q) \leq p^2$. 
 If $Q$ is non-abelian, then 
$\cd(Q) \in \{ \{ 1, p \}, \{ 1, p^2 \}, \{ 1, p, p^2 \} \}$. 
If $\cd(Q) = \{ 1, p \}$, then 
$\mu(G) \leq 2p^3$
(see \cite[Theorem 3.2]{BG} and \cite[Theorem 3]{SAcyclic}).  
If $\cd(Q) \in \{ \{ 1, p^2 \}, \{ 1, p, p^2 \} \}$, then 
$Q \in \Phi_{i}$ where $i\in \{ 5,7,8,10 \}$ (see \cite[\S 4.1]{RJ}). 
Hence $\mu(Q) \in \{ p^2, p^3 \}$ (see \cite[Tables 1--2]{ESAp6}). 
\end{proof}

\begin{lemma} \label{lemma:centerp}
No group 
in $\Phi_{i}$ for \[ i \in \{ 22,24,27,30,
 31, 32,33,35,36, 37, 38, 39, 40,41 \}\] 
is exceptional.
\end{lemma}
\begin{proof}
 We consider the families in turn.  
Let $G$ be a group of order $p^6$.
 \begin{itemize}

\item  $\Phi_{22}$,  $\Phi_{27}$  and $\Phi_{30}$: 
If $G$ is in one of these families, 
then $\mu(G) \in \{ p^3, p^4\}$ 
(see \cite[Tables 3 and 7]{ESAp6}) and $\exp(G) \leq p^2$ (see \cite{Arxiv}). 
Let $Q$ be a quotient of $G$. If $|Q| = p^4$, then $\mu(Q) \leq
 p^3$. Suppose $|Q| = p^5$. 
If $G\in \Phi_{22}$, then $Q = G/Z(G) \cong \phi_{2}(1^5)$; 
 if $G\in \Phi_{27}$, then  $Q = G/Z(G) \cong \phi_{3}(1^5)$;
if  $G\in \Phi_{30}$, then $Q = G/Z(G) \cong \phi_{7}(1^5)$ 
(see \cite[\S 4.1]{RJ}). So $\mu(Q) \leq p^3$. 
Hence $\mu(Q) \leq \mu(G)$ and $G$ has no distinguished quotients. 

 \item $\Phi_{24}$: 
If $G$ is in this family, then 
$\mu(G) = p^3$ (see \cite[Table 3]{ESAp6}) and 
$\exp(G) \leq p^2$ (see \cite{Arxiv}). Let $Q$ be a quotient of $G$. 
If $|Q| = p^4$, then $\mu(Q) \leq p^3$. Suppose $|Q| = p^5$. 
Now $Q = G/Z(G) \cong \phi_{3}(1^5)$ (see \cite[\S 4.1]{RJ}), so 
 $\mu(Q) = p^2 + p$. 
Hence $\mu(Q) \leq \mu(G)$ and 
$G$ has no distinguished quotients. 

 \item  $\Phi_{31}$,  $\Phi_{32}$  and $\Phi_{33}$: 
If $G$ is in one of these families, 
then $\mu(G) \in \{ p^3, p^4 \}$ 
(see \cite[Tables 3 and 7]{ESAp6}) and $\exp(G) \leq p^2$
(see  \cite{Arxiv}). Let $Q$ be a quotient of $G$. 
If $|Q| = p^4$, then $\mu(Q) \leq p^3$. Suppose $|Q| = p^5$. 
Now $Q = G/Z(G) \cong \phi_{4}(1^5)$ (see \cite[\S 4.1]{RJ}), 
so $\mu(Q) = 2p^2$. 
Hence $\mu(Q) \leq \mu(G)$ and 
$G$ has no distinguished quotients. 

 \item $\Phi_{36}$ and $\Phi_{38}$: 
If $G$ is in one of these families, 
then $\mu(G) = p^3$ (see \cite[Table 3]{ESAp6})
 and $\exp(G) \leq p^2$ (see \cite{Arxiv}). Let $Q$ be a quotient of $G$. 
If $|Q| = p^4$, then $\mu(Q) \leq p^3$. Suppose $|Q| = p^5$. 
If $G\in \Phi_{36}$, then $Q = G/Z(G) \cong \phi_{9}(1^5)$; 
if $G\in \Phi_{38}$, then $Q = G/Z(G) \cong \phi_{10}(1^5)$ (see 
\cite[\S 4.1]{RJ}).  So $\mu(Q) \leq p^3$. 
Hence $\mu(Q) \leq \mu(G)$ and 
$G$ has no distinguished quotients. 

\item $\Phi_{37}$ and $\Phi_{39}$: 
If $G$ is in one of these families, 
then $\mu(G) \in \{ p^3, p^4\}$ (see
 \cite[Tables 3 and 7]{ESAp6}) and $\exp(G) \leq p^2$ (see \cite{Arxiv}). 
Let $Q$ be a quotient of $G$. If $|Q| = p^4$, then
 $\mu(Q) \leq p^3$. Suppose $|Q| = p^5$. 
 If $G\in \Phi_{37}$, then $Q = G/Z(G) \cong \phi_{9}(1^5)$;
if $G\in \Phi_{39}$, then $Q = G/Z(G) \cong \phi_{10}(1^5)$ 
(see \cite[\S 4.1]{RJ}). So $\mu(Q) \leq p^3$. 
Hence $\mu(Q) \leq \mu(G)$ and $G$ has no distinguished quotients. 

 \item $\Phi_{40}$ and $\Phi_{41}$: 
If $G$ is in one of these families, 
then $\mu(G) = p^3$ (see \cite[Table 3]{ESAp6})
 and $\exp(G) \leq p^2$ (see \cite{Arxiv}). Let $Q$ be a quotient of $G$. 
If $|Q| = p^4$, then $\mu(Q) \leq p^3$. If $|Q| = p^5$, 
then $Q = G/Z(G) \cong \phi_{9}(1^5)$ (see \cite[\S 4.1]{RJ}). 
So $\mu(Q) = 2p^2$. 
Hence $\mu(Q) \leq \mu(G)$ and $G$ has no distinguished quotients. \qedhere
 \end{itemize}
\end{proof}

 \begin{lemma} \label{lemma:mu(G)p^4notexc}
Let $G$ be a non-abelian group of order $p^6$.
If $\exp(G) \leq p^3$ and $\mu(G) = p^4$, then $G$ is not exceptional.
\end{lemma}
\begin{proof}
If $|Q| = p^4$, then it is not a distinguished quotient of $G$. 
Suppose $|Q| = p^5$. 
If $Q$ is abelian, then $\mu(Q) < p^4$ since $\exp(G) \leq p^3$. 
If $Q$ is non-abelian, then 
$\mu(Q) \leq p^4$ (see \cite[Tables 1--2]{ESAp6}). 
\end{proof}

\begin{lemma}\label{lemma:phi2526}
No group in 
$\Phi_{i}$ for $i \in \{ 25,26,42,43 \}$ is exceptional.
\end{lemma}
\begin{proof}
If $G$ is in one of these families, then $\exp(G) \leq p^3$ (see \cite{Arxiv})
and $\mu(G) = p^4$ (see \cite[Table 7]{ESAp6}). 
Lemma \ref{lemma:mu(G)p^4notexc} shows $G$ is not exceptional.
\end{proof}

\begin{lemma}\label{lemma:phi28}
The only exceptional group in $\Phi_{28}$ is $G_{(28, 2r)}$ where $r=p-1$.
\end{lemma}
\begin{proof}
 If $G \in \Phi_{28}$, then (see \cite[Tables 3 and 7]{ESAp6}) 
 \[ \mu(G) =
 \begin{cases}
 p^3, & \text{if } G = G_{(28,2r)} \text{ where } r = p-1\\
 p^4, & \text{otherwise}.
 \end{cases}
 \]
Let $G$ have a quotient $Q$ of order $p^5$. 
Now $Q = G/Z(G) \cong \phi_{3}(221)b_1$ (see \cite[\S 4.1]{RJ}), 
so $\mu(G/Z(G)) = p^3 + p^2$. 
Therefore, if $G = G_{(28,2r)}$ where $r = p-1$, then 
$\mu(G) < \mu(G/Z(G))$.  If $G$ is any other group in $\Phi_{28}$,
then $Q$ is not distinguished. 
But $\mu(G) = p^4$, thus $G$ has no distinguished quotient of order $p^4$,
and so $G$ is not exceptional.
\end{proof}

We omit the proofs of the next two lemmas, both are 
similar to that of Lemma \ref{lemma:phi28}.

\begin{lemma}\label{lemma:phi29}
The only exceptional group in $\Phi_{29}$ is $G_{(29, 2r)}$ where $r=p-1$.
\end{lemma}

\begin{lemma}\label{lemma:phi34}
The only exceptional group in $\Phi_{34}$ is $G_{(34, 1)}$.
\end{lemma}

\begin{lemma} \label{lemma:QnotExc}
 Let $G$ be a non-abelian group of order $p^6$ satisfying one 
the following: 
 \begin{enumerate}
 \item [$(i)$] $\mu(G) \geq p^4$ and $\exp(G) \leq p^3$, or
 \item [$(ii)$] $\mu(G) = 2p^3$ and $\exp(G) \leq p^2$.
 \end{enumerate}
 Then $G$ is not exceptional.
\end{lemma}

\begin{proof}
Let $1\neq N \lhd G$ and let $Q= G/N$. 
If $Q$ is cyclic, then it is not distinguished.
Suppose $Q$ is not cyclic. If $|Q| \leq p^4$, then $\mu(Q) \leq p^3$
(see \cite[Table 1]{LemieuxExc}). 

Suppose $|Q| = p^5$.  We first prove part (i). 
If $Q$ is abelian, then 
$\mu(Q) \leq p^3+p^2$ since $\exp(Q) \leq p^3$. If $Q$ is non-abelian,
then $\mu(Q) \leq p^4$ (see \cite[Theorem 1.1]{PUp5}). 
Hence, if $\mu(G) \geq p^4$, then $G$ is not exceptional.

We now prove part (ii). Assume that $\mu(G) = 2p^3$ and $\exp(G) \leq p^2$. 
From \cite[Theorem 3.2]{BG} and \cite[Theorem 3]{SAcyclic}, we 
deduce that $\mu(Q) \leq d(Z(Q))p^3$. Thus, if $d(Z(Q)) \leq 2$, 
then $\mu(Q) \leq 2p^3 = \mu(G)$. If $d(Z(Q)) = 3$, 
then $Q \in \Phi_2$ (see \cite[\S 4.1]{RJ}). 
Tables 1 and 2 of \cite{PUp5} show that 
$\mu(Q) \leq p^3+p^2$. Hence 
$G$ is not exceptional.
\end{proof}

\begin{theorem} \label{thm:upperboundexceptional}
 The number of exceptional groups of order $p^6$ is at most
\[ \frac{(33p +467)}{2} + 6\gcd(p-1, 3) + 3\gcd(p-1, 4). \]
\end{theorem}
\begin{proof}
\mbox{}
 \begin{itemize}
 \item By Corollary \ref{cor:phi15}, the $p+3$ groups 
in $\Phi_{15}$ are not exceptional.

 \item By Lemma \ref{lemma:phi21exc}, the 
$(3p^2+4p+5)/2$ groups in $\Phi_{21}$ are not exceptional.

 \item By Lemma \ref{lemma:phi2526}, the $3p+4$ groups in $\Phi_{i}$, 
where $i\in \{ 25,26,42,43 \}$, are not exceptional.

 \item By Lemma \ref{lemma:centerp}, 
the $2p+45+10\gcd(p-1,3)+5\gcd(p-1,4)+2\gcd(p-1,5)$ groups in $\Phi_{i}$, 
where $i \in \{ 22,24,27,30, 31, 32,33,35,36, 37, 38, 39, 40,41 \}$, 
are not exceptional.

 \item By Lemmas \ref{lemma:phi28}--\ref{lemma:phi34}, 
the $2p$ groups in $\Phi_{i}$, where $i\in \{ 28,29,34 \}$, are not exceptional.

 \item By \cite[Tables 1--3 and 5--6]{ESAp6}, 
there are $p+28+2\gcd(p-1,3)+\gcd(p-1,4)$ groups of order $p^6$ in 
$\Phi_{i}$, where $i \in \{ 2,3,4,5,6,7,8,9,10,14 \}$, which have 
minimal degree at least $p^4$. 
By Lemma \ref{lemma:QnotExc}, they are not exceptional. 
We list these groups in Table \ref{t:NExc1}. 
 
 \begin{small}
 
 \begin{longtable}[c]{|l|c|}
 \caption{Groups with minimal degree at least $p^4$ \label{t:NExc1}} \\
 
 \hline
 Group $G$ & $\mu(G)$\\
 \hline
 \endfirsthead
 
 \hline
 \multicolumn{2}{|c|}{Continuation of Table \ref{t:NExc1}}\\
 \hline
 Group $G$ & $\mu(G)$\\
 \hline
 \endhead
 
 \vtop{\hbox{\strut $G_{(2,3)}$ } } & $p^{4}+p$\\
 \hline
 \vtop{\hbox{\strut $G_{(2,4)}$ } } & $p^{4}+p^2$\\
 \hline
 \vtop{\hbox{\strut $G_{(3,1)}$ } } & $p^{4}$\\
 \hline
 \vtop{\hbox{\strut $G_{(3,5r)}$ for $r=1$, $\nu$ } } & $p^{4}+p$\\
 \hline
 \vtop{\hbox{\strut $G_{(3,10r)}$ for $r=1$, $\nu$ } } & $p^{4}+p^2$\\
 \hline
 \vtop{\hbox{\strut $G_{(4,4)}$, $G_{(4,12)}$ } } & $p^{4}+p^2$\\
 \hline
 \vtop{\hbox{\strut $G_{(4,9r)}$ for $r=1$, $\nu$ } } & $p^{4}+p^3$\\
 \hline
 \vtop{\hbox{\strut $G_{(5,1)}$, $G_{(5,2)}$ } } & $p^{4}$\\
 \hline
 \vtop{\hbox{\strut $G_{(6,2r)}$ for $r=1$, $\nu$ } } & $p^{4}+p^2$\\
 \hline
 \vtop{\hbox{\strut $G_{(6,6rs)}$ for $r,s=1$, $\nu$ } } & $p^{4}+p^3$\\
 \hline
 \vtop{\hbox{\strut $G_{(7,1)}$, $G_{(7,2)}$, $G_{(7,4)}$ } } & $p^{4}$\\
 \hline
 \vtop{\hbox{\strut $G_{(7,3r)}$ for $r=1$, $\nu$ } } & $p^{4}$\\
 \hline
 \vtop{\hbox{\strut $G_{(8,1)}$, $G_{(8,3)}$ } } & $p^{4}$\\
 \hline
 \vtop{\hbox{\strut $G_{(8,2r)}$ for $r=0,1,\ldots,p-2$ } } & $p^{4}$\\
 \hline
 \vtop{\hbox{\strut $G_{(9,3r)}$, where $r=1,\omega,\omega^2$ when $p\equiv 1\bmod 3$,} \hbox{\strut and $r=1$ when $p\equiv 2 \bmod 3$ } } & $p^{4}$\\
 \hline
 \vtop{\hbox{\strut $G_{(10,1)}$ } } & $p^{4}$\\
 \hline
 \vtop{\hbox{\strut $G_{(10,2r)}$, where $r=1,\omega,\omega^2,\omega^3$ when $p\equiv 1\bmod 4$,} \hbox{\strut and $r=1,\omega$ when $p\equiv 3 \bmod 4$ } } & $p^{4}$\\
 \hline
 \vtop{\hbox{\strut $G_{(10,3r)}$, where $r=1,\omega,\omega^2$ when $p\equiv 1\bmod 3$,} \hbox{\strut and $r=1$ when $p\equiv 2 \bmod 3$ } } & $p^{4}$\\
 \hline
 \vtop{\hbox{\strut $G_{(14,1)}$, $G_{(14,2)}$ } } & $p^{4}$ \\
 \hline
 \end{longtable}
 \end{small}
 
\item By \cite[Tables 4--6]{ESAp6}, there are  
$(3p^2+23p+56)/2+6\gcd(p-1,3)+2\gcd(p-1,4)$ groups of order $p^6$ in $\Phi_{i}$, 
where $i \in \{ 4,6,11,12,13,16,17,18,19,20,23 \}$, which have 
exponent at most $p^2$ and minimal degree at least $2p^3$. 
By Lemma \ref{lemma:QnotExc}, they are not exceptional.
We list these groups in Table \ref{t:NExc2}. 

 \begin{small}
 
 \begin{longtable}[c]{|l|c|}
 \caption{Groups with exponent at most $p^2$ and 
minimal degree at least~$2p^3$ \label{t:NExc2}} \\
 
 \hline
 Group $G$ & $\mu(G)$\\
 \hline
 \endfirsthead
 
 \hline
 \multicolumn{2}{|c|}{Continuation of Table \ref{t:NExc2}}\\
 \hline
 Group $G$ & $\mu(G)$\\
 \hline
 \endhead
 
 \vtop{\hbox{\strut $G_{(4,3)}$ } } & $2p^{3}$\\
 \hline
 \vtop{\hbox{\strut $G_{(4,23)}$ } } & $2p^{3}+p$\\
 \hline
 \vtop{\hbox{\strut $G_{(4,25r)}$ for $r = \omega, \omega^{3}, \ldots, \omega^{p-2}$} } & $2p^{3}+p$\\
 \hline
 \vtop{\hbox{\strut $G_{(4,32)}$ } } & $2p^{3}+p^2$\\
 \hline
 \vtop{\hbox{\strut $G_{(4,34r)}$ for $r = \omega, \omega^{3}, \ldots, \omega^{p-2}$} } & $2p^{3}+p^2$\\
 \hline
 \vtop{\hbox{\strut $G_{(6,5r)}$ for $r=1$, $\nu$ } } & $2p^{3}$\\
 \hline
 \vtop{\hbox{\strut $G_{(6,18)}$ } } & $2p^{3}+p$\\
 \hline
 \vtop{\hbox{\strut $G_{(6,19r)}$ for $r = \omega, \omega^{3}, \ldots, \omega^{p-2}$ } } & $2p^{3}+p$\\
 \hline
 \vtop{\hbox{\strut $G_{(11,4)}$ } } & $2p^{3}+p^2$\\
 \hline
 \vtop{\hbox{\strut $G_{(11,10r)}$ for $r = 1,2,\ldots,(p-1)/2$ } } & $2p^{3}+p^2$\\
 \hline
 \vtop{\hbox{\strut $G_{(11,14r)}$ for $r = 1$ and $p \equiv 3 \bmod 4$ } } & $2p^{3}+p^2$\\
 \hline
 \vtop{\hbox{\strut $G_{(11,14r)}$ for $r = \nu$ and $p \equiv 1 \bmod 4$ } } & $2p^{3}+p^2$\\
 \hline
 \vtop{\hbox{\strut $G_{(12,6)}$, $G_{(12,9)}$, $G_{(12,10)}$, $G_{(12,13)}$, $G_{(12,15)}$, $G_{(12,16)}$ } } & $2p^{3}$\\
 \hline
 \vtop{\hbox{\strut $G_{(12,12r)}$ for $r = 2,3\ldots,p-1$ } } & $2p^{3}$\\
 \hline
 \vtop{\hbox{\strut $G_{(13,5)}$, $G_{(13,10)}$, $G_{(13,11)}$ } } & $2p^{3}$\\
 \hline
 \vtop{\hbox{\strut $G_{(13,8r)}$ for $r =1, \nu$ } } & $2p^{3}$\\
 \hline
 \vtop{\hbox{\strut $G_{(16,12r)}$ where $r=1,\omega,\omega^2$ when $p\equiv 1\bmod 3$, } \hbox{\strut and $r=1$ when $p\equiv 2 \bmod 3$ } } & $2p^{3}$\\
 \hline
 \vtop{\hbox{\strut $G_{(17,10r)}$ for $r= 1, \nu$ } } & $2p^{3}$\\
 \hline
 \vtop{\hbox{\strut $G_{(17,14r)}$, $G_{(17,15r)}$, $G_{(17,18r)}$, for $r= 1$ or $\nu$ } } & $2p^{3}$\\
 \hline
 \vtop{\hbox{\strut $G_{(17,19rs)}$ where $r= 1$ or $\nu$, $s=1,2,\ldots, p-1$ and $s \neq r^{-1}$} } & $2p^{3}$\\
 \hline
 \vtop{\hbox{\strut $G_{(17,20)}$ } } & $2p^{3}$\\
 \hline
 \vtop{\hbox{\strut $G_{(17,26r)}$ where $r=1,\omega,\omega^2$ when $p\equiv 1\bmod 3$, } \hbox{\strut and $r=1$ when $p\equiv 2 \bmod 3$ } } & $2p^{3}$\\
 \hline
 \vtop{\hbox{\strut $G_{(18,5)}$ } } & $2p^{3}$\\
 \hline
 \vtop{\hbox{\strut $G_{(18,9r)}$ where $r=1,\omega,\omega^2, \omega^3$ when $p\equiv 1\bmod 4$, } \hbox{\strut and $r=1, \nu$ when $p\equiv 3 \bmod 4$ } } & $2p^{3}$\\
 \hline
 \vtop{\hbox{\strut $G_{(18,12r)}$ where $r=1,\omega,\omega^2$ when $p\equiv 1\bmod 3$, } \hbox{\strut 
and $r=1$ when $p\equiv 2 \bmod 3$ } } & $2p^{3}$\\
 \hline
 \vtop{\hbox{\strut $\phi_{19}(2211)d_{0,0,0},$ $\phi_{19}(2211)d_{r,0,t},$ $\phi_{19}(2211)d_{r,s,t},$}  \hbox{\strut $\phi_{19}(2211)f_{r,0},$ $\phi_{19}(2211)f_{r,s},$ $\phi_{19}(2211)g_{r,0,0},$ }\hbox{\strut $\phi_{19}(2211)g_{r,0,t},$ $\phi_{19}(2211)g_{r,s,t}$, } \hbox{where $r$, $s$ and $t$ are 
listed in \cite{Arxiv}} 
 } & $2p^{3}$\\
 \hline
 \vtop{\hbox{\strut $\phi_{19}(2211)h_{r}$, $\phi_{19}(2211)l_{r}$, for $r=1,2,\ldots,p-1$} } & $2p^{3}$\\
 \hline
 \vtop{\hbox{\strut $\phi_{19}(2211)j_{r}$ for $r=1,2,\ldots,(p-1)/2$} } & $2p^{3}$\\
 \hline
 \vtop{ \hbox{\strut $\phi_{19}(2211)k_{r,s}$ for $r= 1,2,\ldots,p-1$ and $s=0,1,\ldots,(p-1)/2$}
\hbox{\strut where $s-r$ and $2r-s$ are 
not divisible by $p$ } } & $2p^{3}$\\
 \hline
 \vtop{\hbox{\strut $\phi_{19}(21^4)f_{r}$, $\phi_{19}(2211)m_{r}$, for $r=1$ or $\nu$} } & $2p^{3}$\\
 \hline
 \vtop{\hbox{\strut $\phi_{19}(21^{4})d_{r,s}$ where $r$ and $s$ are listed in \cite{Arxiv}} 
 } & $2p^{3}$\\
 \hline
 \vtop{\hbox{\strut $\phi_{19}(21^{4})h$ } } & $2p^{3}$\\
 \hline 
 \vtop{\hbox{\strut $G_{(20,5r)}$, $G_{(20,6r)}$ for $r=1$ or $\nu$} } & $2p^{3}$\\
 \hline
 \vtop{\hbox{\strut $G_{(20,8)}$ } } & $2p^{3}$\\
 \hline
 \vtop{\hbox{\strut $G_{(20,12rs)}$ for $r,s = 1, \nu$ } } & $2p^{3}$\\
 \hline
 \vtop{ \hbox{\strut $G_{(20,13rs)}$ for $r,s = 1, \nu$ and $s=1,2,\ldots,p-1$} } & $2p^{3}$\\
 \hline
 \vtop{\hbox{\strut $G_{(20,15r)}$ where $r=1,\omega,\omega^2, \omega^3$ when $p\equiv 1\bmod 4$, } \hbox{\strut and $r=1, \nu$ when $p\equiv 3 \bmod 4$ } } & $2p^{3}$\\
 \hline
 \vtop{\hbox{\strut $G_{(20,16r)}$, $G_{(20,18r)}$, for $r =1,2,\ldots,p-1$ } } & $2p^{3}$\\
 \hline
 \vtop{\hbox{\strut $G_{(20,17r)}$ where $r=1,\omega,\omega^2$ when $p\equiv 1\bmod 3$, } \hbox{\strut and $r=1$ when $p\equiv 2 \bmod 3$ } } & $2p^{3}$\\
 \hline
 \vtop{\hbox{\strut $G_{(23, 8rs)}$ for $s=1$ or $\nu$, and $r=1,\omega,\omega^2$ } \hbox{\strut when $p\equiv 1\bmod 3$, and $r=1$ when $p\equiv 2 \bmod 3$ } } & $2p^{3}$\\
 \hline
 
 \end{longtable}
 \end{small}
 \end{itemize}
We deduce that at least 
\[  3p^2 + \frac{45p+221}{2} + 18\gcd(p-1,3) + 8\gcd(p-1,4)+ 2\gcd(p-1,5) \]
 groups of order $p^6$ are not exceptional. 
The claim follows from (\ref{Total}).
\end{proof}
\cref{thm:main} is an immediate corollary.

\section{A lower bound to the number of exceptional groups of order $p^6$} 
\label{sec:lowerbound}


In Table \ref{t:ExcGroups} we 
identify $(11p+107)/2$ exceptional groups of order $p^6$.  
For each group $G$, we list a generating set for a 
distinguished subgroup $N$.
Group identifiers for $G$ and $G/N$ are those from 
\cite{Arxiv} and \cite{RJ} respectively;
but we identify the distinguished quotient 
of $G_{(4, 38)}, G_{(4, 39)}$ and $G_{(11, 11)}$ as $G_{46}$ 
(see \cite[\S 6.5]{girnat2018}).
The values of $\mu(G)$  and $\mu(G/N)$ are extracted from 
\cite[Tables 2--3]{ESAp6} and  \cite[Tables 1--3]{PUp5} respectively. 

 \begin{scriptsize}

 \begin{longtable}[c]{|l|l|l|l|l|} 
 \caption{Exceptional groups of order $p^6$} 
 \label{t:ExcGroups}\\
 
 \hline
 Group $G$ & \vtop{\hbox{\strut Distinguished} \hbox{\strut subgroup $N$ }} & \vtop{\hbox{\strut Distinguished quotient $Q := G/N$ }} & $\mu(G)$ & $\mu(Q)$ \\
 \hline
 \endfirsthead
 
 \hline
 \multicolumn{5}{|c|}{Continuation of Table \ref{t:ExcGroups}}\\
 \hline
 Group $G$ & \vtop{\hbox{\strut Distinguished} \hbox{\strut subgroup $N$ }} & \vtop{\hbox{\strut Distinguished quotient $Q$ }} & $\mu(G)$ & $\mu(Q)$ \\
 \hline
 \endhead
 
 \hline
 \endfoot
 
 \hline
 \endlastfoot
 
 \multicolumn{5}{|c|}{\footnotesize {\bf Groups in $\Phi_{2}$} }\\
 \hline
 \vtop{\hbox{\strut $G_{(2,7)}$, $G_{(2,8)}$ }} & $\langle \beta_{2}^{-1}\beta_{1}^{p^2} \rangle$ & $\phi_{2}(311)b$ & $p^3+p^2$ & $p^4$ \\
 \hline

 \vtop{\hbox{\strut $G_{(2,20)}$, $G_{(2,22)}$ }} & $\langle \beta_{2}^{-1}\beta_{1}^{p} \rangle$ & $\phi_2(2111)b$ & $2p^2+p$ & $p^3+p$ \\
 \hline

 \vtop{\hbox{\strut $G_{(2,23)}$, $G_{(2,26)}$ }} & $\langle \beta_{2}^{-1}\beta_{1}^{p} \rangle$ & $\phi_2(221)c$ & $3p^2$ & $p^3+p^2$  \\
 \hline

 \multicolumn{5}{|c|}{\footnotesize {\bf Groups in $\Phi_{3}$} }\\
 \hline
 $G_{(3,11)}$, $G_{(3,15)}$ & $\langle \beta_{2}^{-1}\beta_{1}^{p} \rangle$ & $\phi_3(2111)c$ & $2p^2$ & $p^3$ \\
 \hline
 
 $G_{(3,12)}$, $G_{(3,17)}$ & $\langle \beta_{2}^{-1}\beta_{1}^{p} \rangle$ & $\phi_3(311)b_1$ & $p^3+p^2$ & $p^4$ \\
 \hline
 
 $G_{(3,20)}$, $G_{(3,23)}$ & $\langle \beta_{2}^{-1}\beta_{1} \rangle$ & $\phi_3(2111)b_1$ & $2p^2+p$ & $p^3+p$ \\
 \hline 
 
 $G_{(3,25)}$ & $\langle \beta_{1}^{-1}\beta_{3} \rangle$ & $\phi_3(221)b_1$ & $3p^2$ & $p^3+p^2$ \\
 \hline 
 
 \multicolumn{5}{|c|}{\footnotesize {\bf Groups in $\Phi_{4}$}}\\
 \hline
 $G_{(4,2)}$, $G_{(4,10)}$ & $\langle \beta_{2}^{-1}\beta_{1}^{p} \rangle$ & $\phi_2(311)b$ & $p^3+p^2$ & $p^4$ \\
 \hline 
 
 $G_{(4,8r)}$,  for $r=1,2,\ldots,p-1$  & $\langle \beta_{2}^{-1}\beta_{1}^{pr} \rangle$ & $\phi_2(311)b$ & $p^3+p^2$ & $p^4$ \\
 \hline

 $G_{(4,16)}$, $G_{(4, 17)}$ & $\langle \beta_{2}^{-1}\beta_{1} \rangle$ & $\phi_2(2111)b$ & $2p^2+p$ & $p^3+p$ \\
 \hline

 $G_{(4,21r)}$, for $r = \omega, \omega^{2}, \ldots, \omega^{(p-3)/2}$ & $\langle \beta_{2}^{-1}\beta_{1}^{r} \rangle$ & $\phi_2(2111)b$ & $2p^2+p$ & $p^3+p$ \\
 \hline 
 
 $G_{(4,22)}$ & $\langle \beta_{2} \rangle$ & $\phi_2(2111)b$ & $2p^2+p$ & $p^3+p$ \\
 \hline 
 
 $G_{(4,27)}$ & $\langle \beta_{2}^{-1}\beta_{1} \rangle$ & $\phi_2(221)c$ & $3p^2$ & $p^3+p^2$ \\
 \hline 
 
 $G_{(4,30r)}$, for $r = \omega, \omega^{3}, \ldots, \omega^{(p-3)/2}$ & $\langle \beta_{2}^{-1}\beta_{1}^r \rangle$ & $\phi_2(221)c$ & $3p^2$ & $p^3+p^2$ \\
 \hline 
 
 $G_{(4,31)}$ & $\langle \beta_{2} \rangle$ & $\phi_2(221)c$ & $3p^2$ & $p^3+p^2$ \\
 \hline 
 
 $G_{(4,35)}$, $G_{(4, 37)}$ & $\langle \beta_{1}^{-1}\beta_{3} \rangle$ & $\phi_4(2111)c$ & $3p^2$ & $p^3+p^2$ \\
 \hline 
 
 $G_{(4,36)}$ & $\langle \beta_{1}\beta_{2}\beta_{3} \rangle$ & $\phi_4(221)b$ & $3p^2$ & $p^3+p^2$ \\
 \hline

 $G_{(4,38)}$ & $\langle \beta_{1}\beta_{2}^{-1}\beta_{3} \rangle$ & $G_{46}$ & $3p^2$ & $p^3+p^2$ \\
 \hline
 
 $G_{(4,39)}$ & $\langle \beta_{1}\beta_{3} \rangle$ & $G_{46}$ & $3p^2$ & $p^3+p^2$ \\
 \hline
 
 $G_{(4,40)}$, $G_{(4, 41)}$ & $\langle \beta_{1} \rangle$ & $\phi_2(221)c$ & $3p^2$ & $p^3+p^2$ \\
 \hline
 
 \multicolumn{5}{|c|}{\footnotesize {\bf Groups in $\Phi_{6}$} }\\
 \hline
 
 $G_{(6,4)}$ & $\langle \beta_{1}^{p}\beta_{2}^{-1} \rangle$ & $\phi_3(311)b_1$ & $p^3+p^2$ & $p^4$ \\
 \hline
 
 $G_{(6,7r)}$,  for $r =1,2, \ldots, p-1$ & $\langle \beta_{1}^{rp}\beta_{2}^{-1} \rangle$ & $\phi_3(311)b_1$ & $p^3+p^2$ & $p^4$ \\
 \hline 
 
 $G_{(6,9)}$ & $\langle \beta_{1}\beta_{3} \rangle$ & $\phi_6(2111)b_1$ & $3p^2$ & $p^3+p^2$ \\
 \hline 
 
 $G_{(6,10)}$ & $\langle \beta_{1}\beta_{2} \rangle$ & \vtop{\hbox{\strut if $p\equiv 1 
 \bmod 4$ then $\phi_3(2111)b_1$,} {\hbox{\strut if $p\equiv 3 
\bmod 4$ then $\phi_3(2111)b_\nu$ }} } & $2p^2+p$ & $p^3+p$ \\
 \hline
 
 $G_{(6,12r)}$, for $r=1, \nu$ & \vtop{ \hbox{\strut If $p\equiv 3 \bmod 4$ } \hbox{\strut then $\langle \beta_{1}\beta_{3}^{r} \rangle$.} \hbox{\strut If $p\equiv 1 \bmod 4$ } \hbox{\strut and $r=1$ } \hbox{\strut then $\langle \beta_{1}\beta_{3}^{\nu} \rangle$.} \hbox{\strut If $p\equiv 1 \bmod 4$ } \hbox{\strut and $r=\nu$ } \hbox{\strut then $\langle \beta_{1}\beta_{3} \rangle$.}} & $\phi_6(221)d_0$ & $3p^2$ & $2p^3$ \\
 \hline 
 
 $G_{(6,13)}$ & $\langle \beta_{1}\beta_{2}^{-1}\beta_{3} \rangle$ & \vtop{\hbox{\strut if $p\equiv 1 
 \bmod 4$ then $\phi_6(221)c_1$,} {\hbox{\strut if $p\equiv 3 
\bmod 4$ then $\phi_6(221)c_\nu$ }} } & $3p^2$ & $p^3+p^2$ \\
 \hline
 
 $G_{(6,15)}$ & $\langle \beta_{2} \rangle$ & \vtop{\hbox{\strut if $p\equiv 1 
 \bmod 4$ then $\phi_3(2111)b_1$,} {\hbox{\strut if $p\equiv 3 
\bmod 4$ then $\phi_3(2111)b_\nu$ }} } & $2p^2+p$ & $p^3+p$ \\
 \hline 
 
 $G_{(6,16r)}$, for $r=\omega, \omega^2, \ldots, \omega^{(p-3)/2}$ & $\langle \beta_{1}^{-r}\beta_{2} \rangle$ & 
 \vtop{ \hbox{\strut $\phi_3(2111)b_s$, for some $s\in \{ 1, \nu \}$ } }
 & $2p^2+p$ & $p^3+p$ \\
 \hline 
 
 \multicolumn{5}{|c|}{\footnotesize {\bf Groups in $\Phi_{9}$} }\\
 
 \hline
 $G_{(9,8)}$, $G_{(9,10)}$ & $\langle \beta_{1}\beta_{2}^{-1} \rangle$ & $\phi_9(2111)b_1$ & $2p^2$ & $p^3$ \\
 \hline 
 
 \multicolumn{5}{|c|}{\footnotesize {\bf Groups in $\Phi_{11}$}}\\
 \hline
 $G_{(11,3)}$, $G_{(11, 8)}$ & $\langle \alpha_{1} \rangle$ & $\phi_4(2111)c$ & $3p^2$ & $p^3+p^2$ \\
 \hline
 
 $G_{(11,5)}$ & $\langle \alpha_{1}\alpha_{2}^{-1} \rangle$ & $\phi_4(221)f_0$ & $3p^2$ & $2p^3$ \\
 \hline 
 
 $G_{(11,9r)}$, for $r=2,3,\ldots,(p-1)/2$ & $\langle \alpha_{1}\alpha_{2}^{r} \rangle$ & $\phi_4(2111)c$ & $3p^2$ & $p^3+p^2$ \\
 \hline
 
 $G_{(11,11)}$ & $\langle \alpha_{1}\alpha_{3} \rangle$ & $G_{46}$ & $3p^2$ & $p^3+p^2$ \\
 \hline 
 
 $G_{(11,12)}$ & $\langle \alpha_{1}\alpha_{3} \rangle$ & $\phi_4(221)f_s$ for some $s \in \{ 1,2,\ldots,(p-1)/2\}$ & $3p^2$ & $2p^3$ \\
 \hline
 
 $G_{(11,14r)}$ & $\langle \alpha_{1}\alpha_{3} \rangle$ & $\phi_4(221)b$ & $3p^2$ & $p^3+p^2$ \\
 \hline 
 
 $G_{(11,15)}$ & $\langle \alpha_{1}\alpha_{2} \rangle$ & $\phi_4(2111)c$ & $3p^2$ & $p^3+p^2$ \\
 \hline
 
 $G_{(11,16r)}$ & $\langle \alpha_{1}\alpha_{2} \rangle$ & $\phi_4(221)b$ & $3p^2$ & $p^3+p^2$ \\
 \hline
 
 $G_{(11,17r)}$, for $r=1,2,\ldots,(p-1)/2$ & $\langle \alpha_{3} \rangle$ & \vtop{\hbox{\strut $\phi_4(221)f_s$ for some $s\in \{1,2,\ldots,(p-1)/2\}$ } } & $3p^2$ & $2p^3$ \\
 \hline
 
 \multicolumn{5}{|c|}{ \vtop{\hbox{\strut \footnotesize 
{\bf Groups in $\Phi_{12}$}} } }\\
 \hline
 $G_{(12,1)}$ & $\langle \alpha_{1}\alpha_{2} \rangle$ & $\phi_5(1^5)$ & $2p^2$ & $p^3$ \\
 \hline
 
 $G_{(12,2)}$, $G_{(12,8)}$ & $\langle \alpha_{1}\alpha_{2} \rangle$ & $\phi_5(2111)$ & $2p^2$ & $p^3$ \\
 \hline
 
 \multicolumn{5}{|c|}{{\footnotesize {\bf Groups in $\Phi_{16}$}}}\\
 \hline
 $G_{(16,7)}$, $G_{(16,11)}$ & $\langle \alpha_{1}\alpha_{2} \rangle$ & $\phi_3(2111)c$ & $2p^2$ & $p^3$ \\
 \hline 
 
 \multicolumn{5}{|c|}{{\footnotesize {\bf Groups in $\Phi_{17}$} }}\\
 \hline
 $G_{(17,1)}$ & $\langle \alpha_{1}\alpha_{2} \rangle$ & $\phi_7(1^5)$ & $2p^2$ & $p^3$ \\
 \hline
 
 $G_{(17,2)}$ & $\langle \alpha_{1}\alpha_{2} \rangle$ & $\phi_7(2111)a$ & $2p^2$ & $p^3$ \\
 \hline
 
 $G_{(17,6)}$, $G_{(17,16)}$ & $\langle \alpha_{1}\alpha_{2} \rangle$ & \vtop{\hbox{\strut if $p\equiv 1 
 \bmod 4$ then $\phi_7(2111)b_1$,} {\hbox{\strut if $p\equiv 3 
\bmod 4$ then $\phi_7(2111)b_\nu$ }} } & $2p^2$ & $p^3$ \\
 \hline 
 
 $G_{(17,12)}$, $G_{(17,21)}$ & $\langle \alpha_{1}\alpha_{2} \rangle$ & $\phi_7(2111)c$ & $2p^2$ & $p^3$ \\
 \hline 

 \multicolumn{5}{|c|}{ \vtop{\hbox{\strut \footnotesize {\bf Groups in $\Phi_{19}$} } }}\\
 \hline
 
 $\phi_{19}(2211)a$ & $\langle \beta_{1}\beta_{2} \rangle$ & $\phi_7(2111)a$ & $2p^2$ & $p^3$ \\
 \hline
 
 $\phi_{19}(2211)b_{r}$, for $r=1,2, \ldots, (p-1)/2$ & $\langle \beta_{1}\beta_{2}^{r} \rangle$ & \vtop{ \hbox{\strut $\phi_7(2111)b_s$, for some $s\in \{ 1, \nu \}$ } } & $2p^2$ & $p^3$ \\
 \hline
 
 $\phi_{19}(21^4)a$ & $\langle \beta_{1}\beta_{2} \rangle$ & \vtop{\hbox{\strut if $p\equiv 1 
 \bmod 4$ then $\phi_7(2111)b_1$,} {\hbox{\strut if $p\equiv 3 
\bmod 4$ then $\phi_7(2111)b_\nu$ }} } & $2p^2$ & $p^3$ \\
 \hline 
 
 $\phi_{19}(1^6)$ & $\langle \beta_{1}\beta_{2} \rangle$ & $\phi_{7}(1^5)$ & $2p^2$ & $p^3$ \\
 \hline 
 
 \multicolumn{5}{|c|}{{\footnotesize {\bf Groups in $\Phi_{23}$}}}\\
 \hline
 $G_{(23,1)}$ & $\langle \alpha_{1}\alpha_{2} \rangle$ & $\phi_{10}(1^5)$ & $2p^2$ & $p^3$ \\
 \hline
 
 $G_{(23,2)}$ & $\langle \alpha_{1}\alpha_{2} \rangle$ & $\phi_{10}(2111)a_1$ & $2p^2$ & $p^3$ \\
 \hline
 
 $G_{(23,6)}$, $G_{(23, 9r)}$, for $r=0$ & $\langle \alpha_{1}\alpha_{2} \rangle$ & $\phi_{10}(2111)b_1$  & $2p^2$ & $p^3$ \\
 \hline 
 
 \multicolumn{5}{|c|}{{\footnotesize {\bf Groups in $\Phi_{28}$}}}\\
 \hline
 $G_{(28,2r)}$, for $r = p-1$ & $\langle \alpha_{1} \rangle$ & $\phi_3(221)b_1$ 
 & $p^3$ & $p^3+p^2$ \\
 \hline
 
 \multicolumn{5}{|c|}{{\footnotesize {\bf Groups in $\Phi_{29}$}}}\\
 \hline
 $G_{(29,2r)}$, for $r = p-1$ & $\langle \alpha_{1} \rangle$ & $\phi_3(221)b_{\nu}$ 
 & $p^3$ & $p^3+p^2$ \\
 \hline 
 
 \multicolumn{5}{|c|}{{\footnotesize {\bf Groups in $\Phi_{34}$}}}\\
 \hline
 $G_{(34,1)}$ & $\langle \alpha_{1} \rangle$ & $\phi_{4}(221)b$ & $p^3$ & $p^3+p^2$ \\
 \hline
 
 \end{longtable}
     
 \end{scriptsize} 

\begin{table}
 \caption{Number of exceptional groups in some isoclinism families} 
\begin{tabular}{|l|l||l|l||l|l|} 
 \hline
 Isoclinism family & Number & Isoclinism family & Number & Isoclinism family & Number \\
 \hline
 
$\Phi_2$ & 6 & $\Phi_3$ & 7 & $\Phi_4$ & $2p+10$ \\
 \hline
 $\Phi_6$ & $(3p+9)/2$ & $\Phi_9$ & 2 & $\Phi_{11}$ & $(3p+7)/2$ \\
 \hline
 $\Phi_{12}$ & 3 & $\Phi_{16}$ & 2 & $\Phi_{17}$ & 6 \\
 \hline
 $\Phi_{19}$ & $(p+5)/2$ & $\Phi_{23}$ & 4 & $\Phi_{28}$ & 1 \\
 \hline
 $\Phi_{29}$ & 1 & $\Phi_{34}$ & 1 & & \\
 \hline
\end{tabular}
 \label{t:FamilywiseExcGroups}
\end{table}
 

In Table \ref{t:FamilywiseExcGroups} we record the number of 
exceptional groups in those isoclinism families listed in
Table \ref{t:ExcGroups}. 
By summing these values, we deduce that the number of exceptional groups 
of order $p^6$ is at least $(11p+107)/2$.  
Using the machinery of \cite[\S 6]{ESAp6} in {\sc Magma}, we established 
that the number of exceptional groups of order  
$p^6$ for $5 \leq p \leq 13$ is precisely this value.
This provides evidence to support the following.

\begin{conjecture}
The number of exceptional groups of order $p^6$ for primes $p \geq 5$ 
is $(11p+107)/2$.
\end{conjecture}

\section{Access to results}\label{section:access}
The data recorded in Table  \ref{t:ExcGroups} is publicly available in
{\sc Magma} via a GitHub repository \cite{github}.
For each exceptional group in that table, we record 
a presentation and a generating set for a distinguished subgroup. 
The available code can be used to verify our claims for ``small" primes. 

\bibliographystyle{plainurl}
\bibliography{Exceptionalp6}

\end{document}